# A proof of the Gregory–Leibniz series and new series for calculating π

Frank W. K. Firk

The Koerner Center for Emeritus Faculty, Yale University, New Haven CT

**Abstract**

A non-traditional proof of the Gregory-Leibniz series, based on the relationships among the zeta function, Bernoulli coefficients, and the Laurent expansion of the cotangent is given. New series for calculating π are obtained.

## 1. Introduction

There are infinitely many primes in each of the four arithmetic progressions $4n \pm 1$ and $6n \pm 1$ [Dav]. Every twin prime pair $(p, p + 2)$, excluding (3, 5), is of the form $(6n - 1, 6n + 1)$. The Gregory-Leibniz series

$$\pi/4 = 1 - (1/3) + (1/5) - (1/7) + (1/9) - (1/11) + (1/13) - \ldots \quad (1.1)$$

can be rearranged so that products of the form $(4n + 1)(4n - 1)$ occur as denominators

$$\pi/4 = 1 - 2\{(1/15) + (1/63) + (1/143 + (1/255) + \ldots\}$$

$$= 1 - 2\sum_{n=1}^{\infty} 1/((4n)^2 - 1) \quad (1,2)$$

Alternatively, the products of every twin prime pair, excluding (3, 5), appear as denominators in the form:

$$\pi/4 = 1 - (1/3) + (2/(5\cdot7)) + 1/9 - (2/(11\cdot13)) - 1/15 + (2/(17\cdot19)) + \ldots$$

$$= (\sqrt{3}/2)[1 - 2\sum_{n=1}^{\infty} 1/((6n)^2 - 1)]. \text{ (See proof of Theorem 1)} \quad (1.3)$$

Whether or not there are infinitely many twin prime pairs in the Gregory-Leibniz series remains unknown.

In this note, the value of the generalized infinite series

$$S_k = \sum_{n=1}^{\infty} 1/((k \bullet n)^2 - 1),\ k = 3, 4, 5, \ \ldots \quad (1.4)$$

is obtained using relationships among the zeta function $\zeta(2m)$, the Bernoulli coefficients, $B(2m)$, and the Laurent expansion of cot(x) [AS, BS]. A proof of the Gregory-Leibniz series (1.1) is given, and series for π that involve √2, √3, and the sums $S_6$, $S_{12}$, $S_{24}$, and $S_{48}$ are found.

## 2. Theorems

Theorem 2.1

*The infinite sum*

$$S_k = \sum_{n=1}^{\infty} 1/((k \bullet n)^2 - 1),\ k = 3, 4, 5, \ldots$$

$$= [1 - (\pi/k)cot(\pi/k)]/2$$

Proof.

Let

$$S_k = \sum_{n=1}^{\infty} 1/((k \bullet n)^2 - 1),\ k = 3, 4, 5, \ldots \quad (2.1)$$

Factoring out $1/(k \bullet n)^2$, and expanding $(1 - (1/(k \bullet n)^2))^{-1}$ gives

$$S_k = \sum_{n=1}^{\infty} \sum_{m=1}^{\infty} 1/(k \bullet n)^{2m} \quad (2.2)$$

The zeta function $\zeta(2m)$ relates the sum of terms $1/N^{2m}$ to the Bernoulli coefficients, $B_{2m}$ [AS]

$$\zeta(2m) = \sum_{N=1}^{\infty} 1/N^{2m} = (2\pi)^{2m} \mid B_{2m} \mid /(2 \bullet (2m)!),\ m = 1, 2, 3 \ldots \quad (2.3)$$

We can therefore write

$$S_k = \sum_{m=1}^{\infty} 2^{2m} \mid B_{2m} \mid (\pi/k)^{2m} /(2 \bullet (2m)!) \quad (2.4)$$

The Laurent series for $cot(x)$, written in terms of the Bernoulli coefficients $B_{2m}$, is [AS]

$$cot(x) = (1/x) - \sum_{m=1}^{\infty} (2^{2m} \mid B_{2m} \mid x^{2m})/(x \bullet (2m)!) \text{ for } 0 < x < \pi \quad (2.5)$$

Putting $x = \pi/k$, we therefore obtain

$$cot(\pi/k) - (k/\pi) = -\sum_{m=1}^{\infty} (2^{2m} \mid B_{2m} \mid (\pi/k)^{2m} (k/\pi))/(2m)! \quad (2.6)$$

giving

$$S_k = [1 - (\pi/k)cot(\pi/k)]/2,\ k = 3, 4, 5, \ldots \quad \square \quad (2.7)$$

Some values of $S_k$ are

$S_4 = 0.107300918301275845 \ldots$

$S_6 = 0.046550158941445537 \ldots$

$S_{12} = 0.011475691671573332 \ldots$

$S_{24} = 0.002859055853921023 \ldots$

$S_{48} = 0.000714151049012813 \ldots$

The product of every twin prime pair $p(p + 2)$ contributes to the value of the infinite series

$$S_6 = 0.046550158941445537 \ldots$$

(The result

$$\pi \cot(\pi x) = \frac{1}{x} + 2x\sum_{n=1}^{\infty} \frac{1}{x^2 - n^2}$$

is given in Weisstein [Wei] with neither proof nor reference).

Theorem 2.2

$$1 - 2S_4 = \sum_{n=0}^{\infty} (-1)^n/(2n+1), \textit{ the Gregory-Leibniz series.}$$

Proof.

Using (2.7), we have

$$2S_4 = 1 - (\pi/4)cot(\pi/4)$$

$$= 1 - (\pi/4) \qquad (2.8)$$

therefore

$$\pi/4 = 1 - 2\sum_{n=1}^{\infty} 1/((4n)^2 - 1) \qquad (2.9)$$

$$= 1 - 2\{(1/15) + (1/63) + (1/143) + (1/255) + \ldots\}$$

$$= 1 - \{[(1/3) - (1/5)] + [(1/7) - (1/9)] + \ldots\}$$

$$= 1 - (1/3) + (1/5) - (1/7) + (1/9) - \ldots$$

$$= \sum_{n=0}^{\infty} (-1)^n/(2n+1). \quad \square \qquad (2.10)$$

This proof is fundamentally different from the traditional approach that involves expanding the function $1/(1 + x^2)$ as a power series, integrating term-by-term from 0 to 1 (including the remainder after n-terms), letting $n \to \infty$, and equating the result to arctan(x). In such an approach, it is necessary to prove that the method is valid when $x = 1$ in order to use $\arctan(1) = \pi/4$.

## 3. Series for $\pi$

From (2.7), we obtain

$$\pi = k\cdot\tan(\pi/k)[1 - 2S_k],\ k = 3, 4, 5, \ldots \tag{3.1}$$

If the number field is extended to include √2 and √3, the half-angle formulae

$$\sin(\theta/2) = \pm [(1 - \cos(\theta))/2)]^{1/2} \text{ and } \cos(\theta/2) = \pm [(1 + \cos(\theta))/2]^{1/2} \tag{3.2}$$

can be used to give successive values of $\tan(\pi/(6\cdot 2^j))$, $j = 1, 2, 3, \ldots$

where

$$\tan(\pi/6) = 1/\sqrt{3} \tag{3.3}$$

For example,

$$\tan(\pi/12) = 1/(2 + \sqrt{3})$$

$$\tan(\pi/24) = (a - b)/(b - 2),\ a = 2\sqrt{2} \text{ and } b = 1 + \sqrt{3}$$

and

$$\tan(\pi/48) = (2a/(a - b))^{1/2} - (b - 2)/(a - b). \tag{3.4}$$

Using the expression for $\tan(\pi/48)$, we obtain

$$\pi = 48 \times 0.0655435 \ldots (1 - 2S_{48})$$
$$= 3.1460862 \ldots (1 - 2 \times 0.0007142 \ldots)$$
$$= 3.1415926\ldots$$